\documentclass[11pt,leqno]{article}

\usepackage{amssymb,amsmath,amsthm}

\theoremstyle{plain}
\newtheorem{thm}{Theorem}
\newtheorem{cor}[thm]{Corollary}

\newtheorem{lem}[thm]{Lemma}

\theoremstyle{remark}
\newtheorem*{rk}{Remark}

\theoremstyle{definition}

\newcommand{\n}{\noindent}
\newcommand{\ovl}{\overline}
\newcommand{\bb}[1]{\mathbb{#1}}
\newcommand{\cl}[1]{\mathcal{#1}}

\theoremstyle{remark}
\newtheorem{rem}{Remark}
\begin{document}

\title{A similarity degree characterization of nuclear 
$C^*$-algebras}

\author{by\\ \phantom{blank line}\\ Gilles Pisier\footnote{Partially supported by  
NSF grant 
No.~0200690}\\\\ Department of Mathematics\\ 
Texas A\&M University\\ College Station, TX 77843-3368, 
U. S. A.
\\and\\ Universit\'e Paris VI\\ Equipe d'Analyse, Case 186, 75252\\ 
Paris Cedex 05, France\\  {\bf email: pisier@math.tamu.edu}}

\date{}
\maketitle

\abstract{ We show that a $C^*$-algebra $A$ is nuclear iff there is
a number $\alpha<3$ and  a constant $K$
   such that, for any bounded homomorphism
$u\colon\ A \to B(H)$, there is an isomorphism
$\xi\colon\ H\to H$ satisfying $\|\xi^{-1}\|\|\xi\| \le K\|u\|^\alpha$
and such that $ \xi^{-1} u(.) \xi$ is a $*$-homomorphism. In other words, an infinite dimensional
$A$ is nuclear iff its length (in the sense of our previous work on the Kadison similarity problem)
is equal to 2.
}
\bigskip

\n 2000 Math. Subject Classification: Primary 46 L06.
Seconday: 46L07, 46L57

\bigskip

In 1955, Kadison \cite{K} formulated the following conjecture:\  any bounded 
homomorphism $u\colon \ A\to B(H)$, from a $C^*$-algebra into the algebra 
$B(H)$ of all 
bounded operators on a Hilbert space $H$, is similar to a $*$-homomorphism, 
i.e.\ there is an invertible operator $\xi\colon \ H\to H$ such that $x\to \xi 
u(x)\xi^{-1}$ satisfies $\xi u(x^*)\xi^{-1} = (\xi u(x)\xi^{-1})^*$   for all 
$x$ in $A$. This conjecture 
remains unproved, although many partial results are known (see \cite{C1,H4}). In 
particular, by \cite{H4}, we know that $u$ is similar to a $*$-homomorphism iff 
it is completely bounded (c.b.\ in short) in the sense of e.g.\ \cite{Pa} or 
\cite{P4} (to which we refer for background on c.b.\ maps). Moreover, we have
\[
\|u\|_{cb} = \inf\{\|\xi\|\|\xi^{-1}\|\}
\]
where the infimum runs over all invertible $\xi$ such that $\xi 
u(\cdot)\xi^{-1}$ is a $*$-homomorphism. Recall that, by definition, $\|u\|_{cb} 
= \sup\limits_{n\ge 1} \|u_n\|$ where $u_n\colon \ M_n(A)\to M_n(B(H))$ is the 
mapping taking $[a_{ij}]$ to $[u(a_{ij})]$. Thus Kadison's conjecture is 
equivalent to the validity of the implication $\|u\|<\infty \Rightarrow 
\|u\|_{cb}<\infty$. In \cite{P1}, the author proved that if a $C^*$-algebra $A$ 
verifies Kadison's conjecture, then there is a number $\alpha$ for which there 
exists a constant $K$ so that any bounded homomorphism $u\colon \ A\to B(H)$ 
satisfies $\|u\|_{cb}\le K\|u\|^\alpha$. Moreover, the smallest number $\alpha$ 
with this property is an integer denoted by $d(A)$ (and $\alpha=d(A)$ itself 
satisfies the property).

An analogous parameter can be defined for  a discrete group $G$ and it is proved 
in \cite{P1} that $G$ is amenable iff $d(G) \le 2$. The main result of this note 
is the analogous equivalence for $C^*$-algebras:\ a $C^*$-algebra $A$ is nuclear 
(or equivalently amenable, see below) iff $d(A)\le 2$.  In \cite{P1}, we could 
only prove a partial result in this direction.

Let $A,B$ be $C^*$-algebras. Let $\|~~\|_\alpha$ be a $C^*$-norm on their 
algebraic tensor product, denoted by $A\otimes B$; as usual, $A\otimes_\alpha B$ 
then denotes the $C^*$-algebra obtained by completing $A\otimes B$ with respect 
to $\|~~\|_\alpha$. By classical results (see \cite{T}) the set of $C^*$-norms 
admits a minimal and a maximal element denoted respectively by 
$\|\cdot\|_{\min}$ and $\|\cdot\|_{\max}$. Then $A$ is called nuclear if for any 
$B$ we have $A \otimes_{\min} B = A \otimes_{\max} B$, or equivalently 
$\|x\|_{\min} = \|x\|_{\max}$ for any $x$ in $A\otimes B$. We refer the reader 
to \cite{T,KaR} for more information on nuclear $C^*$-algebras. We note in 
particular that by results due to Connes and Haagerup (\cite{C,H}), a 
$C^*$-algebra is nuclear iff it is amenable as a Banach algebra (in 
B.E.~Johnson's sense).

The main result of this note is as follows.

\begin{thm}\label{thm1}
The following properties of $C^*$-algebra $A$ are equivalent:
\begin{itemize}
\item[(i)] $A$ is nuclear.
\item[(ii)] There are  $\alpha<3$ and a constant $K$   such that any bounded 
homomorphism $u\colon \ A\to B(H)$ satisfies $\|u\|_{cb}\le K\|u\|^\alpha$.
\item[(iii)] Same as (ii) with $K=1$ and $\alpha=2$.
\end{itemize}\end{thm}

The implication (i) $\Rightarrow$ (iii) is well known  (see \cite{Bu,C1}).

In the terminology of \cite{P1}, the similarity degree $d(A)$ is the smallest 
$\alpha$ for which the property considered in (ii) above is satisfied. It is 
proved in \cite{P1} that $d(A)$ is always an integer identified as the smallest 
length of a specific kind of factorization for matrices with entries in $A$.

With this terminology, the preceding theorem means that $A$ is nuclear iff 
$d(A)\le 2$. In the infinite dimensional case, $d(A)>1$ hence $A$ is nuclear iff 
$d(A)=2$.

In his   work on derivations (see \cite{C1} and \cite{C2}) Erik 
Christensen isolated the following property $D_k$ for a $C^*$-algebra. Here $k$ 
is any number $\ge 1/2$. A $C^*$-algebra $A$ has property $D_k$ if for any $H$, 
any representation $\pi\colon \ A\to B(H)$, and any $T$ in $B(H)$ the derivation 
$\delta_T\colon \ A\to B(H)$ defined by $\delta_T(a) = \pi(a)T-T\pi(a)$ 
satisfies
\[
\|\delta_T\|_{cb}\le 2k\|\delta_T\|.
\]
With this terminology, Theorem \ref{thm1} implies the following:

\begin{cor}\label{cor11}
Let $A$ be a $C^*$-algebra. The following assertions are equivalent.
\begin{itemize}
\item[(i)] $A$ is nuclear.
\item[(ii)] $A$ satisfies property $D_k$ for some $k<3/2$.
\item[(iii)] $A$ satisfies property $D_1$.
\end{itemize}
\end{cor}

\begin{proof}
Here again the fact that (i) $\Rightarrow$ (iii) is well known (see 
\cite{Bu,C1}). The equivalence between the similarity problem and the derivation 
problem was established by Kirchberg in \cite{Ki}. Refining Kirchberg's 
estimates, the author proved in \cite{P1} (see also \cite[p.~139]{P4}) that 
property $D_k$ implies that the similarity degree $d(A)$ is at most $2k$. Thus 
(ii) $\Rightarrow$ (i) follows from the corresponding implication in Theorem 
\ref{thm1}.
\end{proof}

The main point in Theorem \ref{thm1} is (ii) $\Rightarrow$ (i). In our previous 
work, we could only prove that (ii) implies that $A$ is ``semi-nuclear,'' i.e.\ 
that whenever a representation $\pi\colon \ A\to B(H)$ generates a semifinite 
von~Neumann algebra, the latter is injective. In this note, we show that the 
semifiniteness assumption is not needed. We use the same starting point as in 
\cite{P1},  but we feel the idea of the 
present proof is more transparent than the one in \cite{P1}. In particular, we 
will use the following  result which is part of Th.2.9 in \cite{P2} (obtained independently
in \cite{CS}), but
the latter is inspired by
and  closely related to Haagerup's 
 Th. 2.1 in \cite{H3}.

\begin{thm}\label{thm2}
Let $N\subset B(H)$ be a von Neumann algebra. 
Then $N$ is injective iff there is a constant $C$ such that, for all
$n$, if elements $x_i$ in $N~Ê(i = 1,...,n)$ admit a
decomposition $x_i = \alpha_i + \beta_i$ with $\alpha_i,\beta_i \in B(H)$
such that $\| {\sum \alpha_i^* \alpha_i}\| \leq 1$ and $\| {\sum \beta_i
\beta_i^*} \|\leq 1$, then there is a decomposition $x_i =
 a_i + b_i$ with $a_i, b_i \in N$ such
that $\| {\sum  a_i^*  a_i} \|\leq C^2$ and $\|
{\sum b_i b_i^*} \|\leq C^2$.
\end{thm}

The preceding statement can be viewed as the analogue for von Neumann algebras
of the characterization of amenable discrete groups obtained in \cite{W} (see also \cite{Bo}).

Our main (somewhat) new ingredient is as follows.

\begin{thm}\label{lm4}
Let $M\subset B(H)$ be a von Neumann algebra with a cyclic vector. 
Let $y_1,\ldots, y_n$ in $M'$ be such that for any $x_1,\ldots, x_n$ in $M$ we 
have
\begin{equation}\label{eq110}
\left\|\sum x_iy_i\right\| \le \max\left\{\left\|\sum x^*_ix_i\right\|^{1/2}, 
\left\|\sum x_ix^*_i\right\|^{1/2}\right\}.
\end{equation}
Then there is a decomposition
\[
y_i = a_i+b_i
\]
with $a_i, b_i$ in $M'$ such that
\[
\left\|\sum a_ia_i{}^*\right\|\le 1\quad \text{and}\quad \left\|\sum 
b_i{}^*b_i\right\|\le 1.
\]
\end{thm}

\begin{proof}
We follow a well known kind of argument with roots in \cite{H3}; see also
\cite{Su} and  the 
proof of a theorem due to Kirchberg as presented in \cite[\S 14]{P4} that we 
will follow closely below.  

Recall that the ``row and column" operator spaces
$R_n\subset M_n$ and $C_n\subset M_n$ are defined by:
$$R_n={\rm span}[e_{1i} \mid 1\le i\le n]  \quad C_n={\rm span}[e_{i1} \mid 1\le i\le n]. $$
Let $\Delta_n \subset C_n\oplus R_n$ be the operator 
space spanned by $\delta_i = e_{i1} \oplus e_{1i}$ $(i=1,2,\ldots, n)$. Our 
assumption means that the linear map
\[
v\colon \ \Delta_n\otimes_{\min} M \to  B(H)
\]
defined by
\[
v\left(\sum \delta_i\otimes x_i\right) = \sum x_iy_i
\]
satisfies $\|v\|\le 1$. (Indeed, it is easy to check that the majorant in   
\eqref{eq110} is equal to $\|\sum \delta_i\otimes x_i\|_{\min}$.)

Since $v$ is clearly a two-sided $M$-module map and $M$ has a cyclic vector,  it 
follows by \cite{Sm} (and unpublished work of Haagerup) that $\|v\|_{cb}  = 
\|v\| \le 1$. Let us assume $C_n\oplus R_n\subset B(K)$. Consider now the 
operator system ${\cl S} \subset M_2(B(K) \otimes_{\min} M)$ formed of all 
elements of the form $$\left(\begin{matrix} a&b\\ 
 c^*&d\end{matrix}\right)$$ with $a,d\in 1\otimes M$, and $b,c\in 
\Delta_n\otimes M$. By \cite[Lemma 14.5]{P4}, the mapping
\[
\tilde v\colon \ {\cl S}\to  M_2(B(H))
\]
defined by $$\tilde v \left(\begin{matrix} 1\otimes m_1&b\\ c^*&1\otimes 
m_2\end{matrix}\right) = \left(\begin{matrix} m_1&v(b)\\ 
v(c)^*&m_2\end{matrix}\right)$$ is a unital c.p.\ map
 (or 
  satisfies $\|\tilde v\|_{cb} =1$, which is the same for a unital map).

By Arveson's extension theorem (see
\cite[p. 154]{T3} or  \cite{Pa}), $\tilde v$ admits
a 
unital c.p.\ extension
  $\hat v\colon \ M_2(B(K) \otimes_{\min} M)\to  M_2(B(H))$.
 Then since $\tilde v$ is a 
$*$-homomorphism when restricted to $$\left\{\left(\begin{matrix} m_1&0\\ 
0&m_2\end{matrix}\right)\ \Big| m_i\in M\right\} \simeq M\oplus M,$$ it 
follows by multiplicative domain arguments (see e.g.\ \cite[Lemma 14.3]{P4}) 
that $\hat v$ must be a bimodule map with respect to the natural actions of 
$M\oplus M$. Let $w\colon \ B(K)\otimes_{\min} M\to  B(H)$ be the map 
defined by
\[
w(\beta) = \left[\hat v\left(\begin{matrix} 0&\beta\\ 
0&0\end{matrix}\right)\right]_{12}.
\]
Since $\hat v$ extends $\tilde v$ we have
\[
w\left(\sum \delta_i\otimes x_i\right) = \sum x_iy_i\qquad (x_i\in M).
\]
Moreover, since $\hat v$ is bimodular, for any $z_i$ in $B(K)$
(since
$
z_i\otimes y_i = (z_i\otimes 1) (1\otimes y_i) = (1\otimes y_i) (z_i\otimes 1)
$)
we find
\[
w(z_i\otimes x_i) = w((z_i\otimes 1)) x_i = x_iw(z_i\otimes 1)
\]
and in particular $w(z_i\otimes 1) \in M'$. Thus if we set
\[
a_i = w((e_{i1}\oplus  0)\otimes 1) \quad \text{and}\quad b_i = w((0 \oplus 
e_{1i}) \otimes 1)
\]
then $a_i,b_i\in M'$, $a_i+b_i = w(\delta_i\otimes 1)=x_i$. Finally
\[
\left\|\sum a_ia_i{}^*\right\|^{1/2} = \left\|\sum a_i\otimes e_{1i}\right\| 
\le \|w\|_{cb} \left\|\sum e_{i1} \otimes e_{1i}\right\|\le 1,
\]
and similarly
\[
\left\|\sum b{}^*_ib_i\right\|^{1/2}\le 1.\qquad \qed
\]
\renewcommand{\qed}{}\end{proof}

\begin{rem}\label{rem10} It is easy to see that the preceding result fails without the
cyclicity assumption: Just consider the case $M={\bb C}$ and $M'=B(H)$ with $\dim(H)=\infty$.
\end{rem}
\begin{rem}\label{rem10}
The same proof gives a criterion for a map $u\colon \ E\to M'$ defined on a 
subspace $E\subset A$ of a general $C^*$-algebra $A$ to admit an extension 
$\tilde u\colon \ A\to M'$ with $\|\tilde u\|_{dec}\le 1$. This is essentially 
the same as Kirchberg's \cite[Th.~14.6]{P4}.
\end{rem}

\begin{rem}\label{rem11}
The above Theorem \ref{lm4} may be viewed as an analogue for the operator space $R_n+C_n$ of 
Haagerup's \cite[Lemma 3.5]{H3} devoted to the operator space $\ell^n_1$
equipped with its maximal structure, in the Blecher-Paulsen sense  (see e.g. \cite[\S 3]{P4}).
 Note that while he decomposes into products, we decompose into sums.
\end{rem}

\begin{rem}\label{rem12}
Let $(E_0,E_1)$ be a compatible pair of operator spaces in the sense of \cite[\S 
2.7]{P4}. Then Remark \ref{rem10} gives a sufficient criterion for a map 
$u\colon \ E_0+E_1\to M'$ to admit a decomposition $u = u_0+u_1$ with $u_0\colon 
\ E_0\to M'$ and $u_1\colon \ E_1\to M'$ satisfying $\|u_0\|_{cb}\le 1$ and 
$\|u_1\|_{cb}\le 1$. Assume that $E_0\subset A_0$ and $E_1\subset A_1$, where 
$A_0,A_1$ are $C^*$-algebras, then this criterion actually ensures that there 
are extensions 
\[
\tilde u_0\colon \ A_0\to M'\quad \text{and}\quad \tilde u_1 \colon \ A_1\to 
M'
\]
 with $\|\tilde u_0\|_{dec} \le 1$ and $\|\tilde u_1\|_{dec}\le 1$. In that 
formulation, the converse also holds up to a numerical factor 2.
Note that, in the special case of interest to us, when $E_0=C$ and $E_1=R$,
 then we can take $A_0,A_1$ equal to 
$K(\ell_2)$ (hence nuclear) so that the min and max norms are identical on 
$(A_0\oplus A_1)\otimes M$.
\end{rem}

\n {\bf Notation.} Let $A\subset B(H)$ be any $C^*$-subalgebra.
For any $x_1,\ldots, x_n$ and $y_1,\ldots, y_n$ in $A$, we 
denote
\begin{align}\label{eq101}
&{
\|(x_j)\|_{R\cap C} = \max\left\{\left\|\sum x^*_jx_j\right\|^{1/2}, \left\|\sum 
x_jx^*_j\right\|^{1/2}\right\}}\\
\label{eq102}
&{\|(y_j)\|_{R+C} = \inf\left\{\left\|\sum \alpha^*_j\alpha_j\right\|^{1/2} + 
\left\|\sum
\beta_j\beta^*_j\right\|^{1/2}\right\}},
\end{align}
where the infimum runs over all $\alpha_j,\beta_j$ in $B(H)$ such that $y_j = 
\alpha_j+\beta_j$. Note that, by the injectivity of
$B(H)$, the   definition of $\|(y_j)\|_{R+C}$  does not really
depend on the choice of $H$ or of the embedding $A\subset B(H)$.
The corresponding fact for $   \|(x_j)\|_{R\cap C}$ is obvious.

\begin{cor}\label{cor6}
Let $M\subset B(H)$ be a von Neumann algebra. Then $M$ 
is injective iff there is a constant $C$ such that, for all $n$, all 
$x_1,\ldots, x_n$ in $M$ and $y_1,\ldots, y_n$ in $M'$, we have
\begin{equation}\label{eq103}
\left\|\sum x_iy_i\right\| \le C\|(x_i)\|_{R\cap
C}
\|(y_i)\|_{R+C}.
\end{equation}
\end{cor}

\begin{proof}
If $M$ has a cyclic vector, then this follows immediately from Theorems 
\ref{thm2} 
and \ref{lm4} and the well known fact that
$M'$ is injective iff $M$ is injective
(see \cite[p. 174]{T3}). Now assume that $M$ has a finite cyclic set, i.e.\ there are 
$\xi_1,\ldots, \xi_N$ in $H$ such that $M\xi_1 +\cdots+ M\xi_N$ is dense in $H$. 
Then the vector $(\xi_1,\ldots,\xi_N)$ in $H^N$ is cyclic for $M_N(M) \subset 
M_N(B(H))$. Moreover, it is easy to check that \eqref{eq103} remains true for 
$M_N(M)$ but with $C$ replaced by a constant $C(N)$ (possibly unbounded when $N$ 
grows). Nevertheless, by the first part of the proof it follows that $M_N(M)$ is 
injective and hence, a fortiori, $M$ is injective. 

 In the general case, let $\{\xi_i\mid i\in I\}$ be a dense subset 
of $H$. For any finite subset $J\subset I$, let $H_J$ be the closure
 of \[
\{ \sum\nolimits_{j\in J} a_j(\xi_j)\mid a_j\in
M\}.
\]
Note that $H_J$ is an invariant subspace for
$M$,  so that (since $M$  is self-adjoint)
the orthogonal projection $P_J\colon \ H\to H_J$
belongs to $ M'$.
 Let $\pi_J(a)$ denote the restriction of $a$ to  ${H_J}$.
Then $\pi_J\colon \ M\to B(H_J)$ is a normal
representation, $\pi_J(M)$ admits a finite
cyclic set (namely $\{\xi_i\mid  i\in J\}$), and
it is easy to check that our assumption \eqref{eq103} is still
verified by $\pi_J(M)$ on $H_J$.

Thus, by the first part of the proof,
$\pi_J(M)$ is injective.  This clearly implies
that the von~Neumann algebra $M_J\subset B(H)$
generated by 
$P_J M$ and $I-P_J$ also is injective.
Finally, since $M$ is the weak-$*$ closure of the directed union of the $M_J$'s, 
we conclude that $M$ itself is injective.

 Conversely, if $M$ injective 
then, by Remark \ref{rem4} below, \eqref{eq103} holds with $C=1$.
\end{proof}
\begin{rem}\label{rem4}
Let $M\subset B(H)$ be an injective von~Neumann algebra, so that there is a 
projection $P\colon \ B(H)\to M'$ with
 $\|P\|_{cb}=1$. Then $M$ satisfies \eqref{eq103} with
$C=1$. To see this,  assume
$y_i\in M'$ and 
$\|(y_i)\|_{R+C}<1$, so that $y_i = \alpha_i+\beta_i$ with $\|\Sigma 
\alpha^*_i\alpha_i\|^{1/2} + \|\Sigma \beta_i\beta^*_i\|^{1/2}<1$. Then 
$y_i=a_i+b_i$ with $a_i,b_i\in M'$ satisfying
\[
\left\|\sum a^*_ia_i\right\|^{1/2} + \left\|\sum b_ib^*_i\right\|^{1/2}
\le 
\|P\|_{cb}=1.
\]
Indeed, $a_i=P\alpha_i$ and $b_i = P\beta_i$ clearly verify
this.

\n  Then for any $x_1,\ldots, x_n$ in $M$ we
have by Cauchy-Schwarz
\begin{align*}
\left\|\sum x_ia_i\right\| &\le \left\|\sum x_ix^*_i\right\|^{1/2} \left\|\sum 
a^*_ia_i\right\|^{1/2}\\
\end{align*}
 {  and}
\begin{align*}
\left\|\sum b_ix_i\right\| &\le \left\|\sum x^*_ix_i\right\|^{1/2} \left\|\sum 
b_ib^*_i\right\|^{1/2},
\end{align*}
therefore, since
\begin{align*}
\left\|\sum x_iy_i\right\| &\le \left\|\sum x_ia_i\right\| + \left\|\sum 
b_ix_i\right\|,\\
\end{align*}
 {  we obtain finally}
\begin{align*}
\left\|\sum x_iy_i\right\| &\le \|(x_i)\|_{R\cap C} \|(y_i)\|_{R+C}.
\end{align*}\qed
\end{rem}

We will also use:

\begin{thm}[\cite{P1}]\label{thm5}
A unital operator algebra $A$ satisfies property (ii) in Theorem \ref{thm1} iff 
we have: \ (iv)\ There is a constant $K'$ satisfying the following:\ for any 
linear map $u\colon \ A\to B(H)$ for which there are a Hilbert space $K$,    bounded linear maps 
$v_1,w_1$ from $A$ to $B(K,H)$  and $v_2,w_2$ from $A$ to $B(H,K)$ such that
\begin{equation}\label{eq1}
\forall a,b\in A\qquad u(ab) = v_1(a) v_2(b) + w_1(a)w_2(b)
\end{equation}
we have
\[
\|u\|_{cb} \le K'(\|v_1\| \|v_2\| + \|w_1\| \|w_2\|).
\]
\end{thm}

\begin{rk}
Note that \eqref{eq1} implies that the bilinear map $(a,b)\to  u(ab)$ is c.b.\ 
on  $\max(A)\otimes_h \max(A)$ with c.b.\ norm $\le K'(\|v_1\|\|v_2\| + \|w_1\| 
\|w_2\|)$. Thus, Theorem \ref{thm5} follows from the case $d=2$ of \cite[Th.~4.2]{P1}.
\end{rk}

Another ingredient is the following Lemma which can be derived from \cite{JP} or 
from the more recent paper \cite{PS}.

\begin{lem}\label{lem6}
Let $E$ be a finite dimensional operator space and let $A$ be a $C^*$-algebra. 
Assume that $E$ is a ``maximal'' operator space (equivalently that $E^*$ is a 
minimal one). Then for any c.b.\ map $u\colon \ A\to E$ we have
\[
\forall n\ \forall a_1,\ldots, a_n\in A\qquad \forall\xi_i\in E^*
\]
\begin{equation}\label{eq2}
\left|\sum \langle u(a_j), \xi_j\rangle\right| \le C\|u\|_{cb} \left(\left\|\sum 
a^*_ja_j\right\|^{1/2} + \left\|\sum a_ja^*_j\right\|^{1/2}\right)
 \cdot \sup_{x\in E} \left(\sum |\xi_j(x)|^2\right)^{1/2}
\end{equation}
where $C$ is a numerical constant.
\end{lem}

\begin{proof}
We may apply \cite[Th.\ 1.4]{JP}, arguing as in \cite[Lemma 6.3]{P1} (using \cite[Th.\ 
17.13]{P2} to remove the exactness assumption) this yields \eqref{eq2} with 
$C=2$. Or we may invoke \cite[Th.~0.3]{PS} taking into account \cite[Lemma 
2.3]{PS} (to remove the exactness assumption) and then we again obtain 
\eqref{eq2} with $C=2$.
\end{proof}

For the convenience of the reader, we reproduce here the elementary Lemma 
\ref{lem6.4} from \cite{P1}.

\begin{lem}\label{lem6.4}
Let $(e_i)$ be the canonical basis of the operator
space $\max(\ell_2)$. Let $H$ be any Hilbert space and let $X$ be either
$B({\bb C},H)$ or $B({ H}^*,{\bb C})$, or equivalently let $X$ be
either the column Hilbert space or the row Hilbert space. Then for all
$x_1,\ldots, x_n$ in $X$ we have
$$\left\|\sum\nolimits^n_1 x_i\otimes e_i\right\|_{X\otimes_{\rm min}
\max(\ell_2)} \le \left(\sum \|x_i\|^2\right)^{1/2}.$$
\end{lem}

\begin{proof}
Assume $X = B({\bb C},H)$ or $B(H^*,{\bb C})$. We identify $X$ with
$H$ as a vector space.	Let
$(\delta_m)$ be an orthonormal basis in $H$. Observe that for any finite
sequence $a_m$ in $B(\ell_2)$ we have in both cases
$$\left\|\sum \delta_m \otimes a_m\right\|_{\rm min} \le \left(\sum
\|a_m\|^2\right)^{1/2}.$$
whence we have, for any $x_1,\ldots, x_n$ in $X$,
\begin{align*}
\left\|\sum x_i\otimes e_i\right\| &= \left\|\sum_m \delta_m
\otimes \sum\nolimits_i \langle x_i,\delta_m\rangle e_i\right\|\\
&\le \left(\sum_m \left\|\sum_i \langle x_i,\delta_m\rangle
e_i\right\|^2\right)^{1/2}\\
&= \left(\sum_{m, i} |\langle x_i,\delta_m\rangle|^2\right)^{1/2} =
\left(\sum_i \|x_i\|^2\right)^{1/2}.
\end{align*}
\end{proof}

\begin{proof}[Proof of Theorem \ref{thm1}]
As we already observed,  it suffices to show that (ii) 
implies that $A$ is nuclear. Let $\pi\colon \ A\to B(H)$ be a representation 
and let $M = \pi(A)''$. Using
 Theorem \ref{thm5} and Corollary \ref{cor6}, we will show that (ii) 
implies that $M$ is injective.  
 By the well known results of Choi--Effros and 
Connes (see \cite{CE}), this implies that $A$ is nuclear. Since $\pi(A)\simeq 
A/\text{ker}(\pi)$ is a quotient of $A$, it obviously inherits the property (ii). 
Thus we may as well replace $\pi(A)$ by $A$: we assume $A\subset B(H)$ and let $M=A''$. It suffices to show 
that $M$ is injective.

\n {\bf Claim.} We claim that for any $x_1,\ldots, x_n$ in $M$  and $y_1,\ldots, 
y_n$ in $M'$ we have
\begin{equation}\label{eq10}
\left\| \sum x_j  y_j     \right\| \le 4K'C\|(x_j)\|_{R\cap C} 
\|(y_j)\|_{R+C}.
\end{equation}
\medskip

\n Note:\ It may be worthwhile for the reader to note that $\|(y_j)\|_{R+C}$ is 
(up to a 
factor 2) in operator space duality with $\|(x_j)\|_{R\cap C}$, namely if we set
\[
|||(y_j)||| = \sup\left\{\left\|\sum x_j\otimes y_j\right\|_{\min}\right\}
\]
where the sup runs over all $(x_j)$ in $B(\ell_2)$ such that $\|(x_j)\|_{R\cap 
C}\le 1$, then we have (see e.g.\ \cite{HP})
\[
|||(y_j)||| \le \|(y_j)\|_{R+C} \le 2|||(y_j)|||.
\]

To prove \eqref{eq10} we introduce the operator space $E=\max(\ell^n_2)$, that 
is $n$-dimensional Hilbert space equipped with its ``maximal'' operator space 
structure in the Blecher-Paulsen sense (see \cite[\S 3]{P4}). 
Let us now  fix an $n$-tuple $(y_j)$ in $M'$ such that $\|(y_j)\|_{R+C}<1$.
In addition,  we fix 
$\xi,\eta$ in the unit sphere of $H$.  
Then we define a linear map $u\colon \ M\to E$ as follows:
\[
u(x) = \sum\nolimits_j \langle x  y_j\xi,\eta \rangle e_j
\]
where $e_j$ is the canonical basis of $\ell^n_2$. We will assume that $E\subset 
B(K)$ completely isometrically. The reader may prefer to consider instead of 
$u$, 
the bilinear form $(x,\xi) \to \langle u(x),\xi\rangle$ defined on $M\times E^*$ 
where $E^*$ is now $\ell^n_2$ equipped with its ``minimal'' (or commutative) 
operator space structure obtained by embedding its isometrically into a {\em 
commutative\/} $C^*$-algebra. We will now apply Theorem \ref{thm5} to $u$.  
 
Since we assume $\|(y_j)\|_{R+C}<1$,    we can write
\[
{y_j} = \alpha_j+\beta_j
\]
with $\|\Sigma \alpha^*_j\alpha_j\|<1$ and $\|\Sigma \beta_j\beta^*_j\|<1$. 
Note  that, since  $y_j \in M'$,  for all  $a,b$ in $M$ we have    
\[
{ab}{y_j} = a{y_j}b
\]
and hence
\[
u(ab) = V(a,b) + W(a,b)
\]
where
\begin{align*}
V(a,b) &= \sum \langle  a\alpha_j b\xi,\eta\rangle e_j\\
W(a,b) &= \sum\langle a\beta_j b\xi,\eta\rangle e_j.
\end{align*}
We now claim that we can write for all $a,b$ in $M$
\begin{equation}\label{eq11}
V(a,b) = v_1(a)v_2(b)\quad \text{and}\quad W(a,b) = w_1(a)w_2(b)
\end{equation}
where 
\begin{alignat*}{2}
&v_1\colon \ M\to B(H\otimes K,K), &\quad  &w_1\colon \ M\to B(H\otimes K,K)\\
&v_2\colon \ M\to B(K,H\otimes K), &\quad  &w_2\colon\ M\to B(K, H\otimes K)
\end{alignat*}
are linear maps all with norm $\le 1$.

Indeed, let us set for $h\in H$, $k\in K$ 
\begin{align*}
&v_1(a)(h\otimes k) = \sum\nolimits_j \langle  a\beta_j h, \eta\rangle e_jk\\
&w_1(a)(h\otimes k)= \langle ah, \eta\rangle k\\
&v_2(b)(k) = b\xi \otimes k\\
&w_2(b)(k) = \sum\nolimits_j \alpha_jb\xi\otimes  e_jk.
\end{align*}
Then, it is easy to check \eqref{eq11}. Also, we have trivially
\begin{align*}
\|w_1(a)\| &= \|a^*\eta\| \le \|a\|\\
\|v_2(b)\| &= \|b\xi\| \le \|b\|.
\end{align*}
Moreover, by Lemma \ref{lem6.4} we have
\begin{align*}
\|v_1(a)\|^2 &\le \sum\nolimits_j \|\beta^*_j  a^*\eta\|^2 = \left\langle\sum 
\beta_j\beta^*_j a^* \eta, a^* \eta\right\rangle\\
&\le \|a^*  \eta\|^2 \le \|a\|^2\\
\|w_2(b)\|^2 &\le \sum\nolimits_j\|\alpha_jb\xi\|^2 = \left\langle \sum 
\alpha^*_j\alpha_j b\xi, b\xi\right\rangle\\
&\le \|b\xi\|^2 \le \|b\|^2,
\end{align*}
By  Theorem \ref{thm5}, it follows that
\[
\|u_{|A}\|_{cb} \le 2K'.
\]
Since $u\colon \ M\to B(K)$ is  clearly normal (i.e.\ $\sigma(M,M_*)$ 
continuous) 
and since $A$ is $\sigma(M,M_*)$ dense in $M$, we clearly have
(by the Kaplansky density theorem) 
\[
\|u\|_{cb}= \|u_{|A}\|_{cb}\le 2K'.
\]
Then by Lemma \ref{lem6}, applied with $\xi_j$ biorthogonal to $e_j$, we have 
\[
\forall n\ \forall x_1,\ldots, x_n\in M\qquad
 \left|\langle\sum  x_j  
y_j \xi,\eta\rangle\right|\le 4K'C\|(x_j)\|_{R\cap C}.
\]
Hence,  taking the supremum over
$\xi,\eta$ and using homogeneity,
  we obtain the claimed inequality \eqref{eq10}.
 Then, by Corollary \ref{cor6}, $M$ 
is injective.
\end{proof}
\begin{rk} Since Lemma \ref{lm4} actually holds whenever $A$ is an exact
operator space (with $C$ replaced by twice the exactness constant \cite{JP,PS}),
the proof of Theorem 1 shows that any unital, exact (non selfadjoint) operator algebra
$A\subset B(H)$ with $d(A)\le 2$ in the sense of \cite{P1} satisfies
\eqref{eq103} for some $C$.

\end{rk}

The preceding arguments establish the following result of independent interest.

\begin{thm}\label{thm10}
A $C^*$-algebra $A$ is nuclear iff for any $C^*$-algebra $B$ there is a constant 
$C$ such that, for all $n$, all $x_1,\ldots, x_n$ in $A$ and all $y_1,\ldots, 
y_n$ in $B$ we have
\begin{equation}\label{eq105}
\left\|\sum x_i\otimes y_i\right\|_{\max} \le C\|(x_i)\|_{R\cap C} 
\|(y_i)\|_{R+C}.
\end{equation}
\end{thm}

\begin{proof} 
Let $\pi\colon \ A\to B(H)$ be a representation. Taking $B = \pi(A)'$
(and using the fact that the set of $n$-tuples 
$(x_i)$ in  $A^{**}$ with   $\|(x_i)\|_{R\cap C} \le1$ is the weak-$*$ closure
of its intersection with $A^n$, see e.g. \cite[p. 303]{P4}) we see 
that \eqref{eq105} implies \eqref{eq103}  for $M=\pi(A)''$. Since this holds for any $\pi$, we may 
argue as in the preceding proof (replacing $\pi$ by $\pi_J$) to conclude that 
$\pi(A)''$ is injective, and hence that $A$ is nuclear. Conversely, if $A$ is 
nuclear it is easy to show (see Remark \ref{rem4}) that \eqref{eq105} holds with 
$C=1$.
\end{proof}

\begin{thm}\label{thm11}
A $C^*$-algebra $A$ is nuclear iff for any $C^*$-algebra $B$ there is a constant 
$C$ such that for all $n$, all $x_1,\ldots, x_n$ in $A$ and all $y_1,\ldots, 
y_n$ in $B$ we have
\[
\left\|\sum x_i\otimes y_i\right\|_{\max} \le C\left\|\sum x_i\otimes \bar 
x_i\right\|^{1/2}_{\min} \left\|\sum y_i\otimes \bar y_i\right\|^{1/2}_{\min},
\]
where the min   norms are relative to $A\otimes\ovl A$ and $B\otimes\ovl B$.
\end{thm}

\begin{proof}
It is known (see \cite[(2.12)]{P2}) that $\|\Sigma x_i\otimes \bar x_i\|^{1/2}_{\min} \le 
\|(x_i)\|_{R\cap C}$. Thus, arguing as above, we find that for any 
representation $\pi\colon \ A\to B(H)$ the von~Neumann algebra $M=\pi(A)''$ 
satisfies the following:\ if $y_1,\ldots, y_n$ in $M'$ are such that $\|\Sigma 
y_i\otimes\bar y_i\|_{\min}<1$, then there are $a_i,b_i$ in $M'$ with $y_i = 
a_i+b_i$ such that $\|\Sigma a^*_ia_i\|^{1/2} <C$ and $\|\Sigma 
b_ib^*_i\|^{1/2}<C$. By \cite[Th.~2.9]{P2}, this ensures that $M'$ is 
injective, and hence $A$ is nuclear.
\end{proof}

\begin{rem}\label{rem13}
Note however that by \cite{H1} the inequality
\[
\left\|\sum x_i\otimes \bar x_i\right\|^{1/2}_{\max}\le C\left\|\sum x_i\otimes 
\bar x_i\right\|^{1/2}_{\min}
\]
characterizes the weak expectation property, which is strictly more general than 
nuclearity.
\end{rem}

\begin{rk} It would be nice to know exactly which families of pairs
of operator spaces in duality $(F_n,F_n^*)$ can be used instead of
$F_n=R_n\cap C_n$ or $F_n=OH_n$ to characterize nuclearity (or injectivity)
analogously to the  above Theorems \ref{thm10} and \ref{thm11}
(note that $F_n=R_n$ or $F_n=C_n$
obviously do not work).
\end{rk}

We will say that a function $f\colon \ {\bb N}\to {\bb R}_+$ is ``slowly 
growing'' if, for any $\varepsilon>0$, there is a constant $C_\varepsilon$ such 
that $f(n)\le C_\varepsilon n^\varepsilon$ for all $n\ge 1$.

The rest of the paper is devoted to a technical refinement, based on  
  the following observation: assume that in Theorem 
\ref{thm2} the constant $C$ depends on $n$, i.e. $C=C(n)$ but
that it  is  ``slowly 
growing''. Then $N$ is injective.

 Indeed, as for Theorem 
\ref{thm2},  this observation follows from the same argument as for \cite[Th.~2.9]{P2},
itself based on
  \cite{H3}. Recall Haagerup's characterization of finite injective von Neumann algebras
(\cite[Lemma 2.2]{H3}): $N$ is finite and injective iff
 for any $n$-tuple $(u_i)$ of unitaries  and any 
central projection $p$ in $N$ we have
\begin{equation}\label{eq20} n= \| \sum p u_i \otimes \overline{ p  u_i }\| .\end{equation} 
Actually, for this to hold it suffices
that there exists a  slowly 
growing function $C(n)$ such that for any $n$-tuple $(u_i)$ of unitaries and any 
central projection $p$  in $N$  we have
\begin{equation}\label{eq30} n\le C(n) \| \sum p u_i \otimes \overline{ p  u_i }\|.\end{equation} 
Indeed, if we   set $t= \sum p u_i \otimes \overline{ p  u_i }$ and 
if we apply the preceding inequality to
$(t^*t)^{m}$, take the $m$-th root and let $m$ go to infinity,
then we find that \eqref{eq30} implies  \eqref{eq20} (a similar trick appears in \cite[Lemma 2.2]{H3}).
Given that this is true, the above observation can be deduced, first in the case when $N$ is semifinite, and
then in the general case,
 from the finite case by the same basic reasoning as in \cite{H3}.

 The following theorems are then easy to obtain in the same way as above.

\begin{thm}\label{thm100}
The following properties of a $C^*$-algebra are equivalent.
\begin{itemize}
\item[(i)] $A$ is nuclear.
\item[(ii)] There is a slowly growing function $C\colon \ {\bb N}\to {\bb R}_+$ 
such that for any $n$ and any $C^*$-algebra $B$ we have:
\end{itemize}
\begin{equation}\label{eq201}
\forall (x_i)\in A^n \quad \forall (y_i)\in B^n\qquad \left\|\sum x_i\otimes 
y_i\right\|_{\max}\le C(n) \|(x_i)\|_{R\cap C} \|(y_i)\|_{R+C}.
\end{equation}
\begin{itemize}
\item[(iii)] There is a slowly growing function $C\colon \ {\bb N}\to {\bb R}_+$ 
such that for any $n$ and any $C^*$-algebra $B$, we have:
\end{itemize}
\begin{equation}\label{eq202}
\forall (x_i)\in A^n\quad \forall (y_i)\in B^n\qquad \left\|\sum x_i\otimes 
y_i\right\|_{\max} \le C(n) \left\|\sum x_i\otimes \bar x_i\right\|^{1/2}_{\min} 
\left\|\sum y_i\otimes \bar y_i\right\|^{1/2}_{\min}.
\end{equation}
\end{thm}

\begin{cor}\label{cor111}
A von Neumann algebra $M$ is injective iff there is a slowly growing function
$C\colon \ {\bb N}\to {\bb R}_+$ such that, for any $n$, any mapping $u \colon\ \Delta_n \to M$ admits an extension
$\tilde u \colon\ M_n\oplus M_n  \to M$ such that
$$ \|\tilde u\|_{cb} \le C(n)  \|u\|_{cb}.$$
\end{cor}

\begin{rk} 
Consider a map $u:\ E\to F$ between operator spaces.
Let   $\gamma(u)=\inf\{\|v\|_{cb} \|w\|_{cb}\}$ 
where the infimum runs over all Hilbert spaces $H$ and all possible 
factorizations $u=vw$ of
$u$ through $B(H)$ (here $v:\ B(H)\to F$ and $w:\ E\to B(H)$).
Let $M$ be a von Neumann algebra.
Assume   that there is a constant $C$ so that, for any $n$, any $u:\ R_n\cap C_n\to M$
satisfies $\gamma(u)\le C \|u\|_{cb}$. Then, by the preceding Corollary, $M$ is injective.
Actually, even if $C=C(n)$ depends on $n$, but grows slowly 
when $n\to \infty$, we conclude that $M$ is injective, and hence, a posteriori,  we can
factor through $B(H)$  {\it any} $u$ that takes values in $M$, regardless of its domain.
It seems interesting to investigate which
(sequences of) operator spaces have the property that they
``force" injectivity like $\{R_n\cap C_n\}$. One can show that
$\{ OH_n\}$ has that property too, but obviously not $\{R_n\}$ or $\{C_n\}$,
since these are themselves injective ! 
\end{rk}

\end{document}